\documentclass[11pt]{article}
\usepackage{amssymb,amsmath}
\newcommand{\EE}{\mathcal{E}}

\newcommand{\OO}{\mathcal{O}}

\newcommand{\R}{\mathbb{R}}

\newcommand{\dd}{\,{\rm d}}
\newcommand{\D}{{\rm d}}

\renewcommand{\div}{\mathop{\mathrm{div}}}
\newcommand{\curl}{\mathop{\mathrm{curl}}}

\newcommand{\supp}{\mathop{\mathrm{supp}}}
\renewcommand{\:}{\thinspace :}

\newcommand{\loc}{\mathrm{loc}}

\newtheorem{thm}{Theorem}[section]

\newtheorem{prop}[thm]{Proposition}
\newtheorem{lem}[thm]{Lemma}
\newtheorem{rem}[thm]{Remark}
\newtheorem{cor}[thm]{Corollary}

\newcommand{\QED}{\mbox{}\hfill$\Box$}

\begin{document}

\title{Long-time asymptotics for two-dimensional 
exterior flows with small circulation at infinity}

\author{
\null\\
{\bf Thierry Gallay}\\ 
Institut Fourier\\
UMR CNRS 5582\\
Universit\'e de Grenoble I\\
BP 74\\
38402 Saint-Martin-d'H\`eres, France\\
{\tt Thierry.Gallay@ujf-grenoble.fr}
\and
\\
{\bf Yasunori Maekawa}\\    
Department of Mathematics\\ 
Graduate School of Science\\
Kobe University\\
1-1 Rokkodai, Nada-ku\\          
Kobe 657-8501, Japan\\
{\tt yasunori@math.kobe-u.ac.jp}}

\date{February 17, 2012}

\maketitle

\begin{abstract}
We consider the incompressible Navier-Stokes equations in a
two-dimensional exterior domain $\Omega$, with no-slip boundary
conditions. Our initial data are of the form $u_0 = \alpha \Theta_0 +
v_0$, where $\Theta_0$ is the Oseen vortex with unit circulation at
infinity and $v_0$ is a solenoidal perturbation belonging to 
$L^2(\Omega )^2 \cap L^q(\Omega)^2$ for some $q \in (1,2)$.
If $\alpha \in \R$ is sufficiently small, we show that the
solution behaves asymptotically in time like the self-similar Oseen
vortex with circulation $\alpha$. This is a global stability result,
in the sense that the perturbation $v_0$ can be arbitrarily large, and
our smallness assumption on the circulation $\alpha$ is independent of
the domain $\Omega$.
\end{abstract}

\section{Introduction}\label{sec1}
Let $\Omega \subset \R^2$ be a smooth exterior domain, namely an
unbounded connected open subset of the Euclidean plane with a smooth
compact boundary $\partial \Omega$. We consider the free motion of an
incompressible viscous fluid in $\Omega$, with no-slip boundary
conditions on $\partial\Omega$. The evolution is governed by the
Navier-Stokes equations
\begin{equation}\label{NS}
  \left\{\begin{array}{llll}
  \partial_t u + (u\cdot\nabla)u \,=\, \Delta u - \nabla p~, 
  \quad \div u \,=\, 0~, & \quad & \hbox{for~} x \in \Omega\,, \quad
  &t > 0\,, \\
  u(x,t) \,=\, 0~, & & \hbox{for~} x \in \partial\Omega\,,
  &t > 0\,, \\
  u(x,0) \,=\, u_0(x)~, & & \hbox{for~} x \in \Omega\,,&
  \end{array}\right.
\end{equation}
where $u(x,t) \in \R^2$ denotes the velocity of a fluid particle at
point $x \in \Omega$ and time $t > 0$, and $p(x,t)$ is the pressure in
the fluid at the same point. For simplicity, both the kinematic
viscosity and the density of the fluid have been normalized to $1$.
The initial velocity field $u_0 : \Omega \to \R^2$ is assumed to be 
divergence-free and tangent to the boundary on $\partial\Omega$.

If the initial velocity $u_0$ belongs to the energy space
\[
  L^2_\sigma(\Omega) \,=\, \{u \in L^2(\Omega)^2 \,|\, \div u = 0
  \hbox{ in }\Omega\,,~ u\cdot n = 0 \hbox{ on }\partial\Omega\}~, 
\]
where $n$ denotes here the unit normal on $\partial \Omega$, then it
is known that system \eqref{NS} has a unique global solution $u \in
C^0([0,\infty); L^2_\sigma(\Omega)) \cap C^1((0,\infty); L^2_\sigma(\Omega))
\cap C^0((0,\infty); H^1_0(\Omega)^2\cap H^2(\Omega)^2)$, which satisfies
the energy equality
\[
  \frac12 \|u(\cdot,t)\|_{L^2(\Omega)}^2 + \int_0^t \|\nabla u(\cdot,
  s)\|_{L^2(\Omega)}^2\dd s \,=\, \frac12 \|u_0\|_{L^2(\Omega)}^2~,
  \qquad \hbox{for all } t > 0~.
\]
This result is classical when $\Omega$ is either the whole plane or a
bounded domain of $\R^2$ \cite{Le,KF,FK,CF}. General unbounded
domains, including exterior domains, were considered e.g. in
\cite{KO1}.  It is also known that the kinetic energy $\frac12
\|u(\cdot,t)\|_{L^2(\Omega)}^2$ converges to zero as $t \to \infty$
\cite{Ma,BM,KO1}, and precise decay rates can be obtained under
additional assumptions on the initial data \cite{KO2,CM,BJ}.

In two-dimensional fluid mechanics, however, the assumption that
the velocity field $u$ be square integrable is quite restrictive, 
because it implies (if $u = 0$ on $\partial \Omega$) that the 
associated vorticity field $\omega = \partial_1 u_2 - \partial_2 u_1$
has zero mean over $\Omega$, see \cite[Section~3.1.3]{MB}. In many
important examples, this condition is not satisfied and the kinetic
energy of the flow is therefore infinite. For instance, when 
$\Omega = \R^2$, the Navier-Stokes equations \eqref{NS} have 
a family of explicit self-similar solutions of the form 
$u(x,t) \,=\, \alpha \Theta(x,t)$, $p(x,t) \,=\, \alpha^2 
\Pi(x,t)$, where $\alpha \in \R$ is a parameter and
\begin{equation}\label{Thetadef}
  \Theta(x,t) \,=\, \frac{1}{2\pi}\,\frac{x^\perp}{|x|^2}
   \Bigl(1 - e^{-\frac{|x|^2}{4(1+t)}}\Bigr)~, \qquad
  \nabla \Pi(x,t) \,=\, \frac{x}{|x|^2}|\Theta(x,t)|^2~.
\end{equation}
Here and in the sequel, if $x = (x_1,x_2) \in \R^2$, we denote
$x^\perp = (-x_2,x_1)$ and $|x|^2 = x_1^2 + x_2^2$. The solution
\eqref{Thetadef} is called the {\it Lamb-Oseen vortex} with
circulation $\alpha$. Remark that $|\Theta(x,t)| = \OO(|x|^{-1})$ as
$|x| \to \infty$, so that $\Theta(\cdot,t) \notin L^2(\R^2)^2$, and 
that the circulation at infinity of the vector field $\Theta$ is equal 
to $1$, in the sense that $\oint_{|x|=R} \Theta_1\dd x_1 + \Theta_2\dd x_2 
\to 1$ as $R \to \infty$. The corresponding vorticity distribution
\begin{equation}\label{Xidef}
  \Xi(x,t) \,=\, \partial_1 \Theta_2(x,t) - \partial_2 \Theta_1(x,t) 
  \,=\, \frac{1}{4\pi(1+t)}\,e^{-\frac{|x|^2}{4(1+t)}}~,
\end{equation}
has a constant sign and satisfies $\int_{\R^2} \Xi(x,t)\dd x = 1$ for
all $t \ge 0$. Oseen's vortex plays an important role in the dynamics
of the Navier-Stokes equations in $\R^2$, because it describes the
long-time asymptotics of all solutions whose vorticity distribution is
integrable. This result was first proved by Giga and Kambe for small
solutions \cite{GK}, and subsequently by Carpio for large solutions
with small circulation \cite{Ca}. The general case was finally settled
by Wayne and the first named author \cite{GW2}. It is worth mentioning
that all these results were obtained using the vorticity formulation
of the Navier-Stokes equations.

In the case of an exterior domain $\Omega \subset \R^2$, much less is
known about infinite-energy solutions, mainly because the vorticity
formulation is not convenient anymore due to the boundary
conditions. A general existence result was established by Kozono and
Yamazaki, who proved that system (1) is globally well-posed for
initial data $u_0$ in the weak $L^2$ space $L^{2,\infty}_\sigma(\Omega)$, 
provided that the local singularity of $u_0$ in $L^{2,\infty}$ is
sufficiently small \cite{KY}. In what follows, we consider initial
data of the form
\begin{equation}\label{init}
  u_0 \,=\, \alpha \chi \Theta_0 + v_0~, 
\end{equation}
where $\Theta_0(x) = \Theta(x,0)$ is Oseen's vortex at time $t = 0$,
and $\chi : \R^2 \to [0,1]$ is a smooth, radially symmetric cut-off
function such that $\chi = 0$ on a neighborhood of
$\R^2\setminus\Omega$ and $\chi(x) = 1$ when $|x|$ is sufficiently
large. For any $\alpha \in \R$ and any $v_0 \in L^2_\sigma(\Omega)$,
Theorem~4 in \cite{KY} asserts that the Navier-Stokes equation
\eqref{NS} has a global solution with initial data \eqref{init}, which
is unique in an appropriate class. However, little is known about the
long-time behavior of this solution, and in particular there is no a
priori estimate which guarantees that the $L^{2,\infty}$ norm of $u$
remains bounded for all times.

Very recently, a first result concerning the long-time behavior of 
solutions of \eqref{NS} with initial data of the form \eqref{init} 
was obtained by Iftimie, Karch, and Lacave\:

\begin{thm}\label{IKLthm} {\rm \cite{IKL}} Let $\Omega \subset \R^2$ 
be a smooth exterior domain whose complement $\R^2\setminus\Omega$ is 
a connected set in $\R^2$. For any $v_0 \in L^2_\sigma(\Omega)$, there 
exists a constant $\epsilon = \epsilon(v_0,\Omega) > 0$ such that, 
for all $\alpha \in [-\epsilon,\epsilon]$, the solution of \eqref{NS} 
with initial data \eqref{init} satisfies
\begin{equation}\label{IKLconv}
  \lim_{t \to \infty} t^{\frac12-\frac1p} \|u(\cdot,t) - \alpha 
  \Theta(\cdot,t)\|_{L^p(\Omega)} \,=\, 0~, \qquad 
  \hbox{for all} \quad p \in (2,\infty)~.
\end{equation}
Moreover, there exists $\epsilon_0 = \epsilon_0(\Omega) > 0$ 
such that $\epsilon \ge \epsilon_0$ if $\|v_0\|_{L^2} \le \epsilon_0$.
\end{thm}

Theorem~\ref{IKLthm} shows that solutions of \eqref{NS} which are
finite-energy perturbations of Oseen's vortex $\alpha\Theta_0$ behave
asymptotically in time like the self-similar Oseen vortex $\alpha
\Theta(x,t)$, provided that the circulation at infinity $\alpha$ is
sufficiently small, depending on the size of the initial perturbation.
The conclusion holds in particular when both the circulation $\alpha$
and the finite-energy perturbation $v_0$ are small, so that
Theorem~\ref{IKLthm} extends to exterior domains the result of Giga
and Kambe \cite{GK}. For large solutions, however, the assumption that
$\alpha$ is small depending on $v_0$ is very restrictive. The goal of
the present paper is to prove the following result, which reaches 
a conclusion similar to that of Theorem~\ref{IKLthm} under different 
assumptions on the initial data\:

\begin{thm}\label{main}
Fix $q \in (1,2)$, and let $\mu = 1/q - 1/2$. There exists a constant 
$\epsilon = \epsilon(q) > 0$ such that, for any smooth exterior 
domain $\Omega \subset \R^2$ and for all initial data of the form 
\eqref{init} with $|\alpha| \le \epsilon$ and $v_0 \in L^2_\sigma(\Omega) 
\cap L^q(\Omega)^2$, the solution of the Navier-Stokes equations \eqref{NS} 
satisfies
\begin{equation}\label{conv}
  \|u(\cdot,t) -\alpha\Theta(\cdot,t)\|_{L^2(\Omega)} + 
  t^{1/2} \|\nabla u(\cdot,t) -\alpha\nabla\Theta(\cdot,t)
  \|_{L^2(\Omega)} \,=\, \OO(t^{-\mu})~,
\end{equation}
as $t \to +\infty$.  
\end{thm}

Here, we also suppose that the circulation at infinity is small, and
we assume in addition that the initial perturbation belongs to
$L^2_\sigma(\Omega) \cap L^q(\Omega)^2$ for some $q < 2$. Unlike in
Theorem~\ref{IKLthm}, the limiting case $q = 2$ is not included, and
the proof shows that $\epsilon(q) = \OO(\sqrt{2-q})$ as $q \to 2$.
However, there is absolutely no restriction on the size of the
perturbation $v_0$, hence Theorem~\ref{main} establishes a {\it global
stability property} for the Lamb-Oseen vortices (with small circulation) 
in two-dimensional exterior domains. In this sense, our result can be
considered as a generalization to exterior domains of the work of
Carpio \cite{Ca}, although our proof relies on completely different
ideas. On the other hand, since our perturbations decay faster at
infinity (in space) than those considered by Iftimie, Karch, and
Lacave, we are able to show that the difference $u(x,t) -\alpha
\Theta(x,t)$ converges rapidly to zero, like an inverse power of 
time, as $t \to \infty$. In particular, using \eqref{conv} and 
elementary interpolation, we obtain the estimate
\[
  \sup_{t > 0} t^{\frac1q-\frac1p} \|u(\cdot,t) - \alpha 
  \Theta(\cdot,t)\|_{L^p(\Omega)} \,<\, \infty~, \qquad 
  \hbox{for all} \quad p \in [2,\infty)~,
\]
which improves \eqref{IKLconv} since $q < 2$. 

At this point, it is useful to mention that the assumption that 
$u_0$ can be decomposed as in \eqref{init} for some $\alpha \in 
\R$ and some $v_0 \in L^2_\sigma(\Omega) \cap L^q(\Omega)^2$ is 
automatically satisfied if we suppose that the initial vorticity 
$\omega_0 = \curl u_0$ is sufficiently localized. Indeed, let
us assume for simplicity that $u_0$ vanishes on the boundary 
$\partial \Omega$. For $1 \le p < \infty$, we denote by $\dot 
W^{1,p}_{0,\sigma}(\Omega)$ the completion with respect to 
the norm $u \mapsto \|\nabla u\|_{L^p}$ of the space of all smooth, 
divergence-free vector fields with compact support in $\Omega$. 
Using this notation, we have the following result\:

\begin{prop}\label{vorticity} Fix $q\in (1, 2)$. Assume that 
$u_0$ belongs to $\dot W_{0,\sigma}^{1,p}(\Omega)$ for some $p\in [1,2)$, 
and that the associated vorticity $\omega_0 = {\rm curl ~}u_0$ satisfies 
\begin{equation}\label{omloc}
  \int_{\Omega} (1 + |x|^2)^m |\omega_0 (x)|^2 \dd x \,<\,\infty~, 
\end{equation}
for some $m > 2/q$. If we denote $\alpha = \int_{\Omega} \omega_0 (x) \dd x$, 
then $u_0$ can be decomposed as in \eqref{init} for some 
$v_0 \in L^2_\sigma(\Omega) \cap L^q(\Omega)^2$. In particular, 
if $|\alpha|\leq \epsilon$, the conclusion of Theorem \ref{main} holds.
\end{prop}

For completeness, we give a short proof of Proposition~\ref{vorticity}
in the Appendix. Returning to the discussion of Theorem~\ref{main}, we
emphasize that the smallness condition on the circulation $\alpha$ is
independent of the domain $\Omega$, which can be an arbitrary multiply
connected exterior domain. In fact, the proof will show that the
optimal constant $\epsilon(q)$ is entirely determined by quantities
that appear in the evolution equation for the perturbation of Oseen's
vortex in the whole plane $\R^2$. Note that Oseen vortices are known
to be globally stable for all values of the circulation $\alpha$ when
$\Omega = \R^2$ \cite{GW2}, but in that particular case one can use
the vorticity equation to obtain precise informations on the solutions
of \eqref{NS}. The reader who is not interested in precise convergence 
rates could consider the following variant of Theorem~\ref{main}, where
the condition on the circulation is totally explicit\:

\begin{cor}\label{explicit}
There exists a universal constant $\epsilon_* \ge 4.956$ such that, 
if $|\alpha| < \epsilon_*$ and if $v_0 \in L^2_\sigma(\Omega) \cap 
L^q(\Omega)^2$ for all $q \in (1,2)$, the solution of the Navier-Stokes 
equations \eqref{NS} with initial data \eqref{init} satisfies
$\|u(\cdot,t) -\alpha\Theta(\cdot,t)\|_{L^2(\Omega)} \to 0$ as 
$t \to \infty$.
\end{cor}

The rest of this paper is devoted to the proof of Theorem~\ref{main},
which is quite different from that of Theorem~\ref{IKLthm} in
\cite{IKL}. In the preliminary Section~\ref{sec2}, we collect various
estimates on the truncated Oseen vortex $\chi \Theta$, which can be
verified by direct calculations. In Section~\ref{sec3}, following the
classical approach of Fujita and Kato \cite{FK}, we prove the
existence of a unique global solution of \eqref{NS} for small initial
data of the form \eqref{init}, and we obtain the asymptotics
\eqref{conv} for small solutions. To deal with large solutions, we
derive in Section~\ref{sec4} a ``logarithmic energy estimate'', which
shows that the energy norm of the perturbation $v$ has at most a
logarithmic growth as $t \to \infty$. This is the key new ingredient,
which we use as a substitute for the classical energy inequality which
holds when $\alpha = 0$. Exploiting this estimate and our assumption
that $v_0 \in L^q(\Omega)^2$, we control in Section~\ref{sec5} the
evolution of a fractional primitive of $v$, and we deduce that the
perturbation $v(\cdot,t)$ converges to zero in energy norm, at least
along a sequence of times. Thus we can eventually use the results of
Section~\ref{sec3}, and the conclusion follows.

\section{The truncated Oseen vortex}\label{sec2}

Fix $\rho \ge 1$ large enough so that $\{x \in \R^2 \,|\, |x| \ge \rho\}
\subset \Omega$. Let $\chi(x) = \tilde\chi(x/\rho)$, where $\tilde
\chi  \in C^\infty(\R^2)$ is a radially symmetric cut-off function 
satisfying $\tilde\chi(x) = 0$ when $|x| \le 1$, $\tilde\chi(x) = 1$ 
when $|x| \ge 2$, and $0 \le \tilde\chi(x) \le 1$ for all $x \in \R^2$. 
We define the {\it truncated Oseen vortex} (with unit circulation) 
as follows\:
\begin{equation}\label{uchidef}
  u^\chi(x,t) \,=\,  \chi(x)\Theta(x,t) \,=\, \frac{1}{2\pi}
  \,\frac{x^\perp}{|x|^2} \Bigl(1 - e^{-\frac{|x|^2}{4(1+t)}}\Bigr)\chi(x)~, 
  \qquad x \in \R^2, \quad t \ge 0~.
\end{equation}
Since $\chi$ is radially symmetric and $\supp \chi \subset \{x \in 
\R^2 \,|\, |x| \ge \rho\} \subset \Omega$, it is clear that $u^\chi(x,t)$ 
is a smooth divergence-free vector field which vanishes in a 
neighborhood of $\R^2 \setminus \Omega$. Let $\omega^\chi = 
\partial_1 u^\chi_2 - \partial_2 u^\chi_1$ be the corresponding 
vorticity field, namely
\begin{equation}\label{omegachi}
  \omega^\chi(x,t) \,=\, \chi(x)\Xi(x,t) + \frac{1}{2\pi}\,
  \frac{1}{|x|^2}\Bigl(1 - e^{-\frac{|x|^2}{4(1+t)}}\Bigr)
  x\cdot\nabla\chi(x)~,
\end{equation}
where $\Xi(x,t)$ is defined in \eqref{Xidef}. Since $u^\chi(x,t) =
\Theta(x,t)$ whenever $|x| \ge 2\rho$, the circulation of $u^\chi$ at
infinity is equal to $1$, so that $\int_{\R^2} \omega^\chi \dd x =
1$. Moreover, a direct calculation shows that
\begin{equation}\label{uchirel}
  (u^\chi \cdot \nabla)u^\chi \,=\, \frac12 \nabla|u^\chi|^2 + 
  (u^\chi)^\perp \omega^\chi \,=\, -\frac{x}{|x|^2}|u^\chi|^2~,
\end{equation}
hence there exists a radially symmetric function $p^\chi(x,t)$ such
that $-\nabla p^\chi = (u^\chi \cdot \nabla)u^\chi$. This shows that
$P(u^\chi \cdot \nabla)u^\chi \,=\, 0$, where $P$ denotes the 
Leray-Hopf projection in $\Omega$, namely the orthogonal projection
in $L^2(\Omega)^2$ onto the subspace $L^2_\sigma(\Omega)$. 

The following elementary estimates will be useful\:

\begin{lem}\label{uchiest}\quad \\
{\bf 1.} For any $p \in (2,\infty]$, there exists a constant
$a_p > 0$ such that
\begin{equation}\label{uchi1}
 \|u^\chi(\cdot,t)\|_{L^p(\R^2)} \,\le\, \frac{a_p}{(1+t)^{\frac12 -
 \frac1p}}~, \qquad t \ge 0~.
\end{equation}
{\bf 2.} For any $p \in (1,\infty]$, there exists a constant
$b_p > 0$ such that
\begin{equation}\label{uchi2}
 \|\nabla u^\chi(\cdot,t)\|_{L^p(\R^2)} \,\le\, \frac{b_p}{(1+t)^{1 
 - \frac1p}}~, \qquad t \ge 0~.
\end{equation}
{\bf 3.} For all $t,s \ge 0$, we have
\begin{equation}\label{uchi3}
 \|u^\chi(\cdot,t) - u^\chi(\cdot,s)\|_{L^2(\R^2)}^2 \,\le\, 
 \frac{1}{4\pi}\,\Big|\log \frac{1+t}{1+s}\Big|~. 
\end{equation}
{\bf 4.} There exists a constant $\kappa_1 > 0$ such that, for all 
$t,s \ge 0$,
\begin{equation}\label{uchi4}
 \|\nabla u^\chi(\cdot,t) - \nabla u^\chi(\cdot,s)\|_{L^2(\R^2)}^2 
 \,\le\, \kappa_1\,\Big|\frac{1}{1+t} - \frac{1}{1+s}\Big|~. 
\end{equation}
Moreover all constants $a_p$, $b_p$, and $\kappa_1$ are independent 
of $\rho$, hence of the domain $\Omega$. 
\end{lem}

\noindent{\bf Proof.} By \eqref{uchidef} we have 
\[
  u^\chi(x,t) \,=\, \chi(x)\Theta(x,t) \,=\, \frac{\chi(x)}{\sqrt{1+t}}
  \,\Theta_0\Bigl(\frac{x}{\sqrt{1+t}}\Bigr)~,
\]
where $\Theta_0(x) = \Theta(x,0)$. Since $0 \le \chi \le 1$ and 
$\Theta_0 \in L^p(\R^2)^2$ for all $p > 2$, we find
\[
  \|u^\chi(\cdot,t)\|_{L^p(\R^2)} \,\le\, \frac{1}{\sqrt{1+t}}\, 
  \Bigl\|\Theta_0\Bigl(\frac{\cdot}{\sqrt{1+t}}\Bigr)\Bigr\|_{L^p(\R^2)} 
  \,=\, \frac{\|\Theta_0\|_{L^p(\R^2)}}{(1+t)^{\frac12 -\frac1p}}~, 
  \qquad t \ge 0~.
\]
This proves \eqref{uchi1}. 

Similarly, we have $\partial_i u^\chi = \chi \partial_i\Theta + 
(\partial_i\chi)\Theta$ for $i = 1,2$. As $\partial_i \Theta_0 \in 
L^p(\R^2)^2$ for all $p > 1$, we obtain as before 
\begin{equation}\label{grad1} 
  \|\chi \partial_i\Theta(\cdot,t)\|_{L^p (\R^2)} \,\le\, 
  \frac{1}{1+t}\,\Bigl\|\partial_i\Theta_0\Bigl(\frac{\cdot}{\sqrt{1+t}}
  \Bigr)\Bigr\|_{L^p(\R^2)} \,=\, \frac{\|\partial_i\Theta_0\|_{L^p(\R^2)}}
  {(1+t)^{1-\frac1p}}~,  \qquad t \ge 0~.
\end{equation}
On the other hand, the function $\partial_i\chi$ is supported in 
the annulus $D = \{x \in \R^2 \,|\, \rho \le |x| \le 2\rho\}$, and
satisfies $|\partial_i\chi(x)| \le C\rho^{-1}$ for some $C > 0$ 
independent of $\rho$. Moreover, it follows from \eqref{Thetadef}
that
\[
  |\Theta(x,t)| \,\le\, \frac{1}{2\pi}\min\Bigl(\frac{1}{|x|}\,,\, 
  \frac{|x|}{4(1+t)}\Bigr)~, \qquad x \in \R^2~, \quad t \ge 0~,
\]
hence
\[
  |(\partial_i\chi(x))\Theta(x,t)| \,\le\, C \min\Bigl(\frac{1}{\rho^2}\,,\, 
  \frac{1}{1+t}\Bigr) \mathbf{1}_D(x)~, \qquad x \in \R^2~, 
  \quad t \ge 0~,
\]
where $\mathbf{1}_D$ is the characteristic function of $D$. Taking 
the $L^p$ norm of both sides, we thus obtain
\begin{equation}\label{grad2}
  \|(\partial_i\chi)\Theta(\cdot,t)\|_{L^p(\R^2)} \,\le\, 
  C \rho^{2/p}\min\Bigl(\frac{1}{\rho^2}\,,\, \frac{1}{1+t}\Bigr) 
  \,\le\, \frac{C}{(1+t)^{1-\frac1p}}~, \qquad t \ge 0~.
\end{equation}
Combining \eqref{grad1} and \eqref{grad2}, we arrive at \eqref{uchi2}.

To prove \eqref{uchi3}, we observe that
\begin{align*}
  \|u^\chi(\cdot,t) - u^\chi(\cdot,s)\|_{L^2(\R^2)}^2 \,&\le\,  
  \frac{1}{4\pi^2} \int_{\R^2} \frac{1}{|x|^2}\Bigl(
  e^{-\frac{|x|^2}{4(1+t)}} - e^{-\frac{|x|^2}{4(1+s)}}\Bigr)^2\dd x \\ 
  \,&=\, \frac{1}{2\pi}\log\left\{\frac12 \sqrt{\frac{1+t}{1+s}} + 
  \frac12 \sqrt{\frac{1+s}{1+t}}\right\} \,\le\, \frac{1}{4\pi} 
  \Big|\log \frac{1+t}{1+s}\Big|~,
\end{align*}
for all $t,s \ge 0$. Finally, using \eqref{omegachi}, we find
\[
  \omega^\chi(x,t) -   \omega^\chi(x,s) \,=\, \chi(x)\Bigl(
  \Xi(x,t) - \Xi(x,s)\Bigr) - \frac{x\cdot\nabla\chi(x)}{2\pi|x|^2} 
  \Bigl(e^{-\frac{|x|^2}{4(1+t)}} - e^{-\frac{|x|^2}{4(1+s)}}\Bigr)~.
\]
Thus $\|\nabla u^\chi(\cdot,t) - \nabla u^\chi(\cdot,s)
\|_{L^2(\R^2)}^2 = \|\omega^\chi(\cdot,t) - \omega^\chi(\cdot,s)
\|_{L^2(\R^2)}^2 \le \big (J_1(t,s)^{1/2} + J_2(t,s)^{1/2}\big )^2$, where 
\begin{align*}
  J_1(t,s) \,&=\, \int_{\R^2} \chi(x)^2\Bigl(\Xi(x,t) - \Xi(x,s)\Bigr)^2
  \dd x \,\le\, \int_{\R^2} \Bigl(\Xi(x,t) - \Xi(x,s)\Bigr)^2\dd x \\
  \,&=\, \frac{1}{8\pi}\left\{\frac{1}{1+t} + \frac{1}{1+s}
  - \frac{4}{t+s+2}\right\} \,\le\, \frac{1}{8\pi}\Big|\frac{1}{1+t} - 
  \frac{1}{1+s}\Big|~,
\end{align*}
and
\begin{align*}
  J_2(t,s) \,&=\, \int_{\R^2} \frac{|\nabla\chi(x)|^2}{4\pi^2|x|^2} 
  \Bigl(e^{-\frac{|x|^2}{4(1+t)}} - e^{-\frac{|x|^2}{4(1+s)}}\Bigr)^2\dd x 
  \,\le\, C \rho^{-4} \int_D \Bigl(e^{-\frac{|x|^2}{4(1+t)}} - 
   e^{-\frac{|x|^2}{4(1+s)}}\Bigr)^2\dd x \\ 
   \,&\le\, C \rho^{-2} \sup_{x \in D}\,\Bigl |e^{-\frac{|x|^2}{4(1+t)}} - 
   e^{-\frac{|x|^2}{4(1+s)}}\Bigr| \,\le\, C \Bigl|\frac{1}{1+t} - 
    \frac{1}{1+s}\Bigr|~.
\end{align*}
We thus obtain \eqref{uchi4}, which is the desired estimate. For later 
use, we also observe that $J_2(t,s)$ can be bounded by $C\rho^2(
\frac{1}{1+t} - \frac{1}{1+s})^2$, for some $C > 0$ independent of 
$\rho$. Since $\rho \ge 1$, this gives the alternative estimate
\begin{equation}\label{uchi5}
 \|\nabla u^\chi(\cdot,t) - \nabla u^\chi(\cdot,s)\|_{L^2(\R^2)}^2 
 \,\le\, \frac{1}{8\pi}\,\Big|\frac{1}{1+t} - \frac{1}{1+s}\Big|
 + C\rho^2\,\Big|\frac{1}{1+t} - \frac{1}{1+s}\Big|^{3/2}~,
\end{equation}
which will be used in Section~\ref{sec4}. This concludes the proof 
of Lemma~\ref{uchiest}. \QED

\bigskip
The truncated Oseen vortex is not a solution of the Navier-Stokes
equation, and therefore we need to control the remainder term $R^\chi
= \Delta u^\chi - \partial_t u^\chi = (\Delta\chi)\Theta + 2
(\nabla\chi\cdot \nabla)\Theta$, which has the explicit expression
\begin{equation}\label{Rchidef}
  R^\chi(x,t) \,=\,  \Theta(x,t) \Delta \chi(x) + 2\frac{x\cdot\nabla
  \chi(x)}{|x|^2}\Bigl(x^\perp \Xi(x,t) - \Theta(x,t)\Bigr)~. 
\end{equation}

\begin{lem}\label{Rchiest}
There exists a constant $\kappa_2 > 0$ (independent of $\rho$)
such that, for any $p \in [1,\infty]$,
\begin{equation}\label{Rchi1}
  \|R^\chi(\cdot,t)\|_{L^p(\R^2)} \,\le\, \frac{\kappa_2\, \rho^{\frac{2}{p}
  -1}}{1+t}~,  \qquad t \ge 0~.
\end{equation}
Moreover, for any vector field $u \in H^1_\loc(\R^2)^2$, we have
\begin{equation}\label{Rchi2}
  \Bigl|\int_{\R^2} R^\chi(x,t)\cdot u(x)\dd x\Bigr| 
   \,\le\, \frac{\kappa_2\, \rho}{1+t}\,\|\nabla u\|_{L^2(D)}~, 
  \qquad t \ge 0~,
\end{equation}
where $D = \{x \in \R^2 \,|\, \rho \le |x| \le 2\rho\}$. 
\end{lem}

\noindent{\bf Proof.} 
It is clear from \eqref{Rchidef} that $|R^\chi(x,t)| \le C\rho^{-1}
(1+t)^{-1}\mathbf{1}_D(x)$ for all $x \in \R^2$ and all $t \ge 0$,  
and \eqref{Rchi1} follows immediately. Moreover, we have 
$R^\chi(x,t) = x^\perp Q^\chi(x,t)$ for some radially symmetric 
scalar function $Q(x,t)$, hence $R^\chi(\cdot,t)$ has zero mean
over the annulus $D$. If $u \in H^1_\loc(\R^2)^2$ and if we denote by
$\bar u$ the average of $u$ over $D$, the Poincar\'e-Wirtinger 
inequality implies
\[
  \Bigl|\int_{\R^2} R^\chi(x,t)\cdot u(x)\dd x\Bigr| \,=\,
  \Bigl|\int_D R^\chi(x,t)\cdot (u(x) - \bar u)\dd x\Bigr| 
  \,\le\, C \rho \|R^\chi(\cdot,t)\|_{L^2(\R^2)} \|\nabla u\|_{L^2(D)}~, 
\]
and using \eqref{Rchi1} with $p = 2$ we obtain \eqref{Rchi2}. 
\QED

\section{Asymptotic behavior of small solutions}\label{sec3}

Given $\alpha \in \R$, we consider solutions of \eqref{NS} of the 
form
\begin{equation}\label{udec}
  u(x,t) \,=\, \alpha u^\chi(x,t) + v(x,t)~, \qquad
  p(x,t) \,=\, \alpha^2 p^\chi(x,t) + q(x,t)~,
\end{equation}
where $u^\chi(x,t)$ is the truncated Oseen vortex \eqref{uchidef}
and $p^\chi$ is the associated pressure. The perturbation $v(x,t)$ 
satisfies the no-slip boundary condition and the equation
\begin{equation}\label{veq}
  \partial_t v + \alpha (u^\chi\cdot\nabla)v + \alpha (v\cdot\nabla)u^\chi
  + (v\cdot\nabla)v \,=\, \Delta v + \alpha R^\chi - \nabla q~, 
  \qquad \div v \,=\, 0~,
\end{equation}
where $R^\chi$ is given by \eqref{Rchidef}. If we apply the Leray-Hopf 
projection $P$ and use the fact that $PR^\chi = R^\chi$, we obtain 
the equivalent system
\begin{equation}\label{veq2}
  \partial_t v + \alpha P \Bigl((u^\chi\cdot\nabla)v + (v\cdot\nabla)
  u^\chi\Bigr) + P(v\cdot\nabla)v \,=\, -A v + \alpha R^\chi~,
\end{equation}
where $A = -P\Delta$ is the Stokes operator, which is selfadjoint 
and nonnegative in $L^2_\sigma(\Omega)$ with domain $D(A) = 
L^2_\sigma(\Omega ) \cap H^1_{0}(\Omega)^2 \cap H^2(\Omega )^2$.

In this section, we fix some initial time $t_0 \ge 0$ and prove 
the existence of global solutions to \eqref{veq2} with small 
initial data $v_0 = v(\cdot,t_0)$ in the energy space. The 
integral equation associated with \eqref{veq2} is 
\begin{align}\nonumber
   v(t) \,=\, S(t-t_0)v_0 + \int_{t_0}^t S(t-s)& \Bigl\{
   \alpha R^\chi(s) - P(v(s)\cdot\nabla)v(s) \\ \label{veq3}
   &-\alpha P\Bigl((u^\chi(s)\cdot\nabla)v(s) + (v(s)\cdot\nabla)
   u^\chi(s)\Bigr)\Bigr\}\dd s~, 
\end{align}
where $v(t) \equiv v(\cdot,t)$ and $S(t) = \exp(-tA)$ is the Stokes 
semigroup. For $p \in (1,\infty)$, we denote by $L^p_\sigma(\Omega)$
the closure in $L^p(\Omega)^2$ of the set of all smooth divergence-free
vector fields with compact support in $\Omega$. We then have the following 
standard estimates\:

\begin{prop}\label{Stokes} The Stokes operator $-A$ generates an analytic 
semigroup of contractions in $L^2_\sigma(\Omega)$. Moreover, for each 
$t > 0$ the operator $S(t) = \exp(-tA)$ extends to a bounded linear 
operator from $L^q_\sigma(\Omega)$ into $L^2_\sigma(\Omega)$ for $1 < q 
\le 2$, and there exists a constant $C = C(q) > 0$ (independent 
of $\Omega$) such that
\begin{equation}\label{Stokesbd}
  t^{\frac1q-\frac12} \|S(t)v_0\|_{L^2(\Omega)} + t^{\frac1q} 
  \|\nabla S(t)v_0\|_{L^2(\Omega)} \,\le\, C \|v_0\|_{L^q}~, \qquad t > 0~,
\end{equation}
for all $v_0 \in L^q_\sigma(\Omega)$. In particular, we can take 
$C = 2$ in \eqref{Stokesbd} if $q = 2$.
\end{prop}

Since $A$ is selfadjoint and nonnegative, it is clear that
$\{S(t)\}_{t\ge 0}$ is an analytic semigroup of contractions in
$L^2_\sigma(\Omega)$, and that both terms in the left-hand side of
\eqref{Stokesbd} are bounded by $\|v_0\|_{L^2}$. On the other hand, 
general $L^q-L^p$ estimates for $S(t)$ were established in 
\cite{BV,DS1,DS2,KY,MS}, but the corresponding constants depend 
a priori on the domain $\Omega$. The fact that \eqref{Stokesbd} 
holds with $C$ independent of $\Omega$ was essentially observed 
in \cite{BM,KO2}. For the reader's convenience, we repeat the proof 
of \eqref{Stokesbd} in Section \ref{sec5} below.

The main result of this section is\:
\begin{prop}\label{small} There exist positive constants $K_0$, $\delta$, 
$V_\Omega$, and $T_\Omega$ such that, if $t_0 \ge T_\Omega$, if $|\alpha|\le 
\delta$, and if $\|v_0\|_{L^2(\Omega)} \le V_\Omega$, then the perturbation 
equation \eqref{veq2} has a unique global solution $v \in C^0([t_0,
\infty); L^2_\sigma(\Omega))$ such that
\begin{equation}\label{vbound}
  \sup_{t \ge t_0} \|v(t)\|_{L^2(\Omega)} +  \sup_{t > t_0}(t-t_0)^{\frac12} 
  \|\nabla v(t)\|_{L^2(\Omega)} \,\le\, 4 \|v_0\|_{L^2(\Omega)} + 
  K_0\,\rho^{\frac12} |\alpha| (1+t_0)^{-\frac14}~. 
\end{equation}
Here $K_0$ and $\delta$ are independent of $\Omega$. In addition,  
if there exists $\mu \in (0,1/2)$ such that 
\begin{equation}\label{moredecay}
  M \,:=\, \sup_{\tau > 0}\tau^\mu \|S(\tau)v_0\|_{L^2(\Omega)} + 
  \sup_{\tau > 0}\tau^{\mu + \frac12} \|\nabla S(\tau)v_0\|_{L^2(\Omega)} 
  \,<\, \infty~, 
\end{equation}
then 
\begin{equation}\label{vdecay}
  \sup_{t > t_0}(t-t_0)^\mu \|v(t)\|_{L^2(\Omega)} +  \sup_{t > t_0}
  (t-t_0)^{\mu + \frac12} \|\nabla v(t)\|_{L^2(\Omega)} \,\le\, 2 M + 
  C_\Omega |\alpha|~,
\end{equation}
for some $C_\Omega > 0$ depending on $\Omega$. 
\end{prop}

\noindent{\bf Proof.} We follow the classical approach of Fujita
and Kato \cite{FK}. Given $t_0 \ge 0$, we introduce the Banach space 
$X = \{v \in C^0([t_0,\infty); L^2_\sigma(\Omega)) \cap C^0((t_0,\infty); 
H^1_0(\Omega)^2) \,|\, \|v\|_X < \infty\}$, equipped with the norm
\[
  \|v\|_X \,=\, \sup_{t \ge t_0} \|v(t)\|_{L^2} +  \sup_{t > t_0}
  (t-t_0)^{\frac12}\|\nabla v(t)\|_{L^2}~.
\]
If $v_0 \in L^2_\sigma(\Omega)$, we denote $\bar v(t) = S(t-t_0)v_0$
for $t \ge t_0$. In view of \eqref{Stokesbd}, we have $\bar v \in X$ 
and $\|\bar v\|_X \le 2\|v_0\|_{L^2}$. On the other hand, given any 
$v \in X$ we denote, for $t \ge t_0$, 
\[
  (F v)(t) \,=\, \int_{t_0}^t S(t-s) \big ( \alpha R^\chi(s)  + \alpha G_1^v(s) 
  + G_2^v (s) \big ) \dd s \, = \, \alpha F_0 (t) + \alpha (F_1 v) (t) 
  + (F_2 v) (t)~,
\]
where $G_1^v(s) = - P(u^\chi(s)\cdot\nabla)v(s) -P(v(s)\cdot\nabla)u^\chi(s)$ 
and $G_2^v (s) = -P(v(s)\cdot\nabla)v(s)$. We shall show that $F$ maps $X$ 
into $X$, and that there exist positive constants $C_1, C_2, C_{3,\Omega}$
(independent of $t_0$) such that
\begin{align}\label{Fest1}
  \|Fv\|_X \,&\le\, C_1 \rho^\frac{1}{2} |\alpha| (1+t_0 )^{-\frac14} + |\alpha| 
  C_2  \|v\|_X + C_{3,\Omega} \|v\|_X^2~, \\ \label{Fest2}
  \|Fv - F\tilde v\|_X \,&\le\,  |\alpha|  C_2 \|v -\tilde v\|_X +
  C_{3,\Omega} (\|v\|_X + \|\tilde v\|_X) \|v- \tilde v\|_X~, 
\end{align}
for all $v,\tilde v \in X$. 

To prove \eqref{Fest1}, we estimate separately the contributions 
of $F_0$, $F_1$, and $F_2$. First, using \eqref{Stokesbd} with $q = 4/3$, 
we obtain for $t > t_0$\:
\begin{equation}\label{Fest3.0}
  \|F_0(t)\|_{L^2} + (t{-}t_0)^{\frac12}\|\nabla F_0(t)\|_{L^2} \,\le\, 
  C \int_{t_0}^t \left(\frac{1}{(t{-}s)^{\frac14}} + \frac{(t{-}t_0)^{\frac12}}{
  (t{-}s)^{\frac34}}\right)\| R^\chi (s)\|_{L^{\frac43}}\dd s~,
\end{equation}
and from Lemma~\ref{Rchiest} we know that $\|R^\chi(s)\|_{L^{4/3}} \le
C\rho^{1/2} (1+s)^{-1}$ for all $s \ge 0$.  It follows that $\|F_0\|_X
\le C_1 \rho^{1/2}(1+t_0 )^{-1/4}$ for some $C_1 > 0$ independent of
$t_0$ and $\Omega$. In a similar way, we find
\begin{equation}\label{Fest3.2}
  \|(F_2 v)(t)\|_{L^2} + (t{-}t_0)^{\frac12}\|\nabla (F_2 v)(t)\|_{L^2} 
  \,\le\, C  \int_{t_0}^t \left(\frac{1}{(t{-}s)^{\frac14}} + 
  \frac{(t{-}t_0)^{\frac12}}{(t{-}s)^{\frac34}}\right)\| G_2 ^v 
  (s)\|_{L^{\frac43}}\dd s~.
\end{equation}
Using the fact that the Leray-Hopf projection is a bounded operator
in $L^{4/3}(\Omega)^2$, whose norm depends a priori on $\Omega$, 
we estimate
\[
 \|G_2^v(s)\|_{L^{\frac43}} \,\le\, C_\Omega \|v(s)\|_{L^4}
 \|\nabla v(s)\|_{L^2} \,\le\, C_\Omega \|v(s)\|_{L^2}^{\frac12}\|\nabla v(s)
 \|_{L^2}^{\frac32} \,\le\, \frac{C_\Omega \|v\|_X^2}{(s-t_0)^{\frac34}}~,
\]
for all $s > t_0$. It follows that $\|F_2 v\|_X \le C_{3,\Omega}
\|v\|_X^2$, where $C_{3,\Omega} > 0$ is independent of $t_0$. Finally,
to bound $F_1$, we proceed in a slightly different way in order
to obtain a constant $C_2$ that does not depend on $\Omega$. Observing
that $G_1^v (s) = A^{1/2} A^{-1/2} P \div (u^\chi\otimes v + v\otimes u^\chi)
(s)$, and that $\|A^{1/2}v\|_{L^2} = \|\nabla v\|_{L^2}$ for all $v \in
L^2_\sigma(\Omega)) \cap H^1_0(\Omega)^2$, we can use \eqref{Stokesbd} 
with $q = 2$ to obtain
\begin{equation}\label{Fest3.1}
  \|(F_1 v)(t)\|_{L^2}  \,\le\, \int_{t_0}^t (t{-}s)^{-\frac12} \|A^{-1/2} 
  P \div (u^\chi\otimes v  + v\otimes u^\chi)(s)\|_{L^{2}}\dd s~.
\end{equation}
Similarly, the quantity $(t-t_0)^\frac12 \|\nabla (F_1 v)(t)\|_{L^2}$
can be bounded by
\begin{equation}\label{Fest3.1'}
  \int_{t_0}^{\frac{t+t_0}{2}} \frac{(t{-}t_0)^\frac12}{t{-}s} \|A^{-1/2} P 
  \div (u^\chi\otimes v  + v\otimes u^\chi)(s)\|_{L^{2}}\dd s + 
  \int_{\frac{t+t_0}{2}}^t \frac{(t{-}t_0)^\frac12}{(t{-}s)^{\frac12}} \|G_1^v(s)
  \|_{L^{2}}\dd s~.
\end{equation}
Since $A^{-1/2}P\div$ defines a bounded operator from $L^2(\Omega)^4$ 
into $L^2_\sigma (\Omega )$ whose norm is less than or equal to $1$ 
(see \cite[Lemma III-2-6-1]{Sohr}), we have from \eqref{uchi1}
\[
  \|A^{-1/2} P \div (u^\chi\otimes v  + v\otimes u^\chi)(s)\|_{L^{2}} 
  \,\le\, 2 \|u^\chi (s)v(s)\|_{L^2} \,\le\, 2 a_\infty (1+s)^{-\frac12} 
  \|v\|_X~.
\]
Similarly, using \eqref{uchi1} and \eqref{uchi2} we find
\[
  \|G_1^v (s)\|_{L^{2}} \,\le\, \|u^\chi(s) \nabla v(s)\|_{L^2} + 
  \|v(s) \nabla u^\chi(s)\|_{L^2} \,\le\, \frac{a_\infty \|v\|_X}{
  (1+s)^{\frac12}(s-t_0)^{\frac12}} + \frac{b_\infty \|v\|_X}{1+s}~.
\]
Inserting these estimates into \eqref{Fest3.1} and \eqref{Fest3.1'}, 
we obtain $\|F_1 v\|_X \le C_2 \|v\|_X$ for some $C_2 > 0$ independent
of $t_0$ and $\Omega$. Since $Fv = \alpha F_0 + \alpha F_1v + F_2v$, 
this concludes the proof of \eqref{Fest1}, and the Lipschitz bound 
\eqref{Fest2} is established in exactly the same way.

Now let $B_r = \{v \in X\,|\, \|v\|_X \le r\}$, where $r > 0$ is small
enough so that $4 r C_{3,\Omega} \le 1$. If we assume that $4 |\alpha|
C_2 \le 1$, $8\|v_0\|_{L^2} \le r$, and $4C_1\rho^{1/2}|\alpha|(1+t_0)^{-1/4}
\le r$, the estimates above imply that the map $v \mapsto \bar v -
Fv$ leaves the closed ball $B_r$ invariant and is a strict contraction
in $B_r$. By construction, the unique fixed point of that map in $B_r$
is the desired solution of \eqref{veq3}. This proves the existence
part of Proposition~\ref{small} with
\[
  K_0 \,=\, 2 C_1~, \quad \delta \,=\, \frac{1}{4C_2}~, \quad
  V_{\Omega} \,=\, \frac{1}{32 C_{3,\Omega}}~, \quad 
  T_\Omega \,=\,\Bigl(\frac{4 C_1 C_{3,\Omega}\rho^{\frac12}}{C_2}\Bigr)^4~.
\]

In a second step, we assume that \eqref{moredecay} holds for some
$\mu \in (0,1/2)$. Given any $T > t_0$, we denote
\[
  \EE_T \,=\, \sup_{t_0 \le t \le T}(t-t_0)^\mu \|v(t)\|_{L^2} +  
  \sup_{t_0 < t \le T}(t-t_0)^{\mu + \frac12} \|\nabla v(t)\|_{L^2}~, 
\]
where $v$ is the solution of \eqref{veq3} constructed in the 
previous step. Our goal is to show that $\EE_T$ is uniformly bounded by 
a constant which does not depend on $T$. From \eqref{veq3}, we 
know that
\begin{equation}\label{EET}
  \EE_T \,\le\, M + \sup_{t_0 \le t \le T}(t-t_0)^\mu \|(Fv)(t)\|_{L^2} +  
  \sup_{t_0 < t \le T}(t-t_0)^{\mu + \frac12} \|\nabla (Fv)(t)\|_{L^2}~, 
\end{equation}
where $M$ is defined in \eqref{moredecay}. To estimate the last two 
terms, we proceed as above. Let $p \in (1,2)$ be such that $1/p > 
\mu + 1/2$, and define $q \in (2,\infty)$ by the relation $1/q = 
1/p - 1/2$. As in \eqref{Fest3.0} and \eqref{Fest3.2}, we have
\begin{align*}
  & (t{-}t_0)^\mu \|F_0(t)\|_{L^2} + (t{-}t_0)^{\mu + \frac12}\|\nabla F_0(t)
  \|_{L^2} \,\le\, C \int_{t_0}^t \left(\frac{ (t{-}t_0)^\mu}{(t{-}s)^{\frac1q}} 
  + \frac{(t{-}t_0)^{\mu + \frac12}}{(t{-}s)^{\frac1p}}\right)\| R^\chi(s)
  \|_{L^{p}}\dd s~, \\
  & (t{-}t_0)^\mu \|( F_2 v ) (t)\|_{L^2} + (t{-}t_0)^{\mu + \frac12}\|\nabla 
  (F_2v)(t)\|_{L^2} \le C \!\int_{t_0}^t \left(\frac{(t{-}t_0)^\mu}{
  (t{-}s)^{\frac1q}} + \frac{(t{-}t_0)^{\mu + \frac12}}{(t{-}s)^{\frac1p}}
  \right)\| G_2 ^v (s)\|_{L^{p}}\dd s~,
\end{align*}
for $t \in (t_0,T]$. Moreover $\|R^\chi(s)\|_{L^p} \le C \rho^{\frac2p-1}
(1+s)^{-1}$ and
\[
 \|P(v(s)\cdot\nabla)v(s)\|_{L^p} \,\le\, C_\Omega  \|v(s)\|_{L^q}
 \|\nabla v(s)\|_{L^2} \,\le\, C_\Omega  \|v(s)\|_{L^2}^{\frac2q}\|\nabla v(s)
 \|_{L^2}^{2-\frac2q} \,\le\, \frac{C_\Omega \|v\|_X \EE_T}{
 (s-t_0)^{\mu+1-\frac1q}}~,
\]
for all $s \in (t_0,T]$. The term involving $F_1v$ is estimated as in 
\eqref{Fest3.1} and \eqref{Fest3.1'}, and we find
\begin{align*}
  (t-t_0)^\mu \|(F_1 v)(t)\|_{L^2}  \,&\le\,  C \int_{t_0}^t 
  \frac{(t-t_0)^\mu \EE_T}{(t{-}s)^{\frac12} (1+s)^\frac12 (s-t_0)^\mu}  
  \dd s~,\\
  (t-t_0)^{\mu + \frac12} \|\nabla (F_1 v)(t)\|_{L^2}  \,&\le\, C 
  \int_{t_0}^\frac{t+t_0}{2} \frac{(t-t_0)^{\mu + \frac12} \EE_T }{(t-s) 
  (1+s)^{\frac12} (s-t_0)^{{\mu}}} \dd s \\
  &+ C \int_{\frac{t+t_0}{2}}^t \frac{(t-t_0)^{\mu + \frac12}}{(t{-}s)^{\frac12}} 
  \Bigl(\frac{\EE_T }{(1+s)^\frac12 (s-t_0)^{{\mu +\frac12}}} 
  + \frac{\EE_T}{(1+s) (s-t_0)^\mu }\Bigr)\dd s~.
\end{align*}
If we insert these estimates into \eqref{EET}, we obtain after 
elementary calculations
\begin{equation}\label{EET2}
   \EE_T \,\le\, M + \tilde C_1 \rho^{\frac{2}{p}-1} |\alpha| 
   (1+t_0)^{-\frac1p +\mu + \frac12} + \tilde C_2 |\alpha| \EE_T
  + \tilde C_{3,\Omega} \|v\|_X \EE_T~,
\end{equation}
for some positive constants $\tilde C_1, \tilde C_2, \tilde C_{3,\Omega}$ 
independent of $T$ and $t_0$. Now, taking $\delta$ and $V_\Omega$ smaller 
and $T_\Omega$ larger if needed, we can ensure that $\tilde C_2 
|\alpha| + \tilde C_{3,\Omega} \|v\|_X \le 1/2$. Then \eqref{EET2} implies 
that
\[ 
  \EE_T \,\le\, 2M + 2 \frac{\tilde C_1 \rho^{\frac2p-1} |\alpha|}
  {(1+t_0)^{\frac1p -\mu -\frac12 }}~,
\]
for all $T > t_0$, and \eqref{vdecay} follows. This concludes the
proof. \QED

\begin{rem}\label{locrem}
{\rm The proof of Proposition~\ref{small} can be modified in a classical
way \cite{FK,Br} to yield the following local existence result. For 
any $\alpha \in \R$, any $t_0 \ge 0$, and any $v_0 \in L^2_\sigma
(\Omega)$, there exists $T = T(\alpha,v_0,\Omega) > 0$ such that
Eq.~\eqref{veq2} has a unique solution $v \in C^0([t_0,t_0+T];
L^2_\sigma(\Omega)) \cap C^0((t_0,t_0+T]; H^1_0(\Omega)^2)$ satisfying
$v(t_0) = v_0$; moreover, any upper bound on $|\alpha| + 
\|v_0\|_{H^1}$ gives a lower bound on the local existence time $T$.
In our formulation of Proposition~\ref{small}, smallness conditions
were imposed on $\alpha$ and $v_0$ to ensure global existence, and 
the assumption on the intial time $t_0$ guarantees that the 
smallness condition on $\alpha$ is independent of the domain $\Omega$.} 
\end{rem}

\section{A logarithmic energy estimate}\label{sec4}

In this section, we establish our key estimate for large solutions of
\eqref{veq2} in the energy space. Fix $\alpha \in \R$, $v_0 \in 
L^2_\sigma(\Omega)$, and let $v \in C^0([0,T];L^2_\sigma(\Omega)) \cap 
C^0((0,T];H^1_0(\Omega)^2)$ be a solution of \eqref{veq2} with initial data 
$v(0) = v_0$, see Remark~\ref{locrem}. We first derive a crude bound on 
$v$ using a classical energy estimate. Multiplying both sides of 
\eqref{veq2} by $v$ and integrating by parts over $\Omega$, we find
\begin{equation}\label{venergy}
  \frac12 \frac{\D}{\D t} \|v(t)\|_{L^2}^2 + \|\nabla v(t)\|_{L^2}^2 
  \,=\, \alpha \langle v(t), R^\chi(t)\rangle - \alpha \langle v(t), 
  (v(t)\cdot\nabla) u^\chi(t)\rangle~,
\end{equation}
where $\langle \cdot,\cdot \rangle$ denotes the usual scalar product 
in $L^2_\sigma(\Omega)$, so that $\|\cdot\|_{L^2} = \langle \cdot\,,\cdot
\rangle^{1/2}$. Using \eqref{Rchi2}, we easily obtain 
\[
  |\alpha \langle v(t),R^\chi(t)\rangle| \,\le\, \frac{\kappa_2\, \rho 
  |\alpha|}{1+t} \|\nabla v(t)\|_{L^2} \,\le\, \frac{\eta}{2} 
  \|\nabla v(t)\|_{L^2}^2 + \frac{\kappa_2^2 \rho^2 \alpha^2}{2\eta(1+t)^2}~,
\]
for any $\eta \in (0,1]$. Moreover, applying \eqref{uchi2} with $p = 
\infty$, we see that 
\[
  |\langle v(t), (v(t)\cdot\nabla) u^\chi(t)\rangle| \,\le\, 
  \frac{b_\infty}{1+t}\,\|v(t)\|_{L^2}^2~.
\]
We thus obtain the energy inequality
\[
  \frac{\D}{\D t} \|v(t)\|_{L^2}^2 + (2 - \eta  ) \|\nabla v(t)\|_{L^2}^2
  \,\le\,  \frac{2b_\infty |\alpha|}{1+t}\|v(t)\|_{L^2}^2 + 
  \frac{\kappa_2^2 \rho^2 \alpha^2}{\eta(1+t)^2}~, \qquad 0 < t \le T~.
\]
Using Gronwall's lemma, we deduce that
\begin{equation}\label{vpoly}
  \|v(t)\|_{L^2}^2 +  (2 - \eta ) \int_{t_0}^t \|\nabla v(s)\|_{L^2}^2\dd s 
  \,\le\, \Bigl(\frac{1+t}{1+t_0}\Bigr)^{2 b_\infty |\alpha|}
  \Bigl(\|v(t_0)\|_{L^2}^2 + \frac{\kappa_2^2 \rho^2 \alpha^2}{\eta(1+t_0)} 
  \Bigr)~,
\end{equation}
for $0 \le t_0 < t \le T$. 

We shall see that estimate \eqref{vpoly} is pessimistic for large
times, but it already implies that the solutions of \eqref{veq2} in
the energy space $L^2_\sigma(\Omega)$ are global. Indeed,
\eqref{vpoly} shows that the norm $\|v(t)\|_{L^2}$ grows at most
polynomially in time, and it is then straightforward to establish a similar
result for $\|\nabla v(t)\|_{L^2}$. In particular, the $H^1$ norm of
$v(t)$ cannot blow up in finite time, and using Remark~\ref{locrem} 
we conclude that all solutions of \eqref{veq2} in $L^2_\sigma(\Omega)$ 
are global.

The aim of this section is to establish the following ``logarithmic 
energy estimate'', which improves \eqref{vpoly} for large times. 

\begin{prop}\label{logbound} There exists a constant $K_1 > 0$ (independent
of $\Omega$) such that, for any $\alpha \in \R$ and any $v_0 \in 
L^2_\sigma(\Omega)$, the solution of \eqref{veq2} with initial data $v_0$ 
satisfies, for all $t \ge 1$,
\begin{equation}\label{vlog}
  \|v(t)\|_{L^2(\Omega)}^2 + \int_0^t \|\nabla v(s)\|_{L^2(\Omega)}^2\dd s 
  \,\le\, K_1 \Bigl(\|v_0\|_{L^2(\Omega)}^2 + \alpha^2 \log(1+t) + 
  D_{\alpha,\rho}\Bigr)~,
\end{equation}
where $D_{\alpha,\rho} = \alpha^2\log(1+|\alpha|) + \alpha^2 \rho^2$.  
\end{prop}

\noindent{\bf Proof.} As in \eqref{vpoly}, we introduce here a 
parameter $\eta \in (0,1]$, which will be used in Section~\ref{sec5}
below to specify the optimal smallness condition on the circulation
$\alpha$ and prove Corollary~\ref{explicit}. The reader who is not
interested in optimal constants should set $\eta = 1$ everywhere. 

Given any $\tau \ge 0$, we denote
\begin{equation}\label{tildevdef}
  \tilde v(x,t) \,=\, u(x,t) - \alpha u^\chi(x,t+\tau) \,=\, 
  v(x,t) + \alpha \Bigl(u^\chi(x,t)- u^\chi(x,t+\tau)\Bigr)~, 
\end{equation}
for all $x \in \Omega$ and all $t > 0$. Then $\tilde v$ satisfies
\eqref{veq2} where $u^\chi(x,t)$ and $R^\chi(x,t)$ are replaced by
$u^\chi(x,t+\tau)$ and $R^\chi(x,t+\tau)$, respectively. Proceeding
exactly as above, we thus obtain the following energy estimate\:
\begin{equation}\label{tildevpoly}
  \|\tilde v(t)\|_{L^2}^2 +  (2 - \eta ) \int_0^t  \|\nabla \tilde v(s)
  \|_{L^2}^2\dd s \,\le\, \Bigl(\frac{1+t+\tau}{1+\tau}\Bigr)^{2 b_\infty 
  |\alpha|} \Bigl(\|\tilde v(0)\|_{L^2}^2 + \frac{\kappa_2^2 \rho^2 
  \alpha^2}{\eta (1+\tau) } \Bigr)~,
\end{equation}
for all $t > 0$. Now, we fix $t \ge 1$ and choose $\tau = Nt-1$, where
\[
  N \,=\, N_{\alpha,\eta} \,=\, \max\Bigl(1\,,\,\frac{2b_\infty|\alpha|}{
  \log(1+\eta)}\Bigr)~.
\]
This choice implies that
\[
  \Bigl(\frac{1+t+\tau}{1+\tau}\Bigr)^{2 b_\infty |\alpha|} \,=\, 
  \Bigl( 1+\frac{1}{N}\Bigr)^{2 b_\infty |\alpha|} \,\le\, 1+\eta~.   
\]
On the other hand, using \eqref{uchi3}, \eqref{tildevdef}, we find
\begin{align*}
  \|v(t)\|_{L^2}^2  \,&\le\, (1{+}\eta) \|\tilde v(t)\|_{L^2}^2 + 
  \frac{1{+}\eta}{\eta} \alpha^2\|u^\chi(t) - u^\chi(t{+}\tau)\|_{L^2}^2 
  \,\le\, (1{+}\eta)\|\tilde v(t)\|_{L^2}^2 + \frac{\alpha^2}{2\pi\eta } 
  \log (N{+}1)~, \\
  \|\tilde v(0)\|_{L^2}^2 \,&\le\, \frac{1{+}\eta}{\eta}\|v_0\|_{L^2}^2 + 
  (1{+}\eta)\alpha^2 \|u^\chi(0) - u^\chi(\tau)\|_{L^2}^2
  \,\le\,\frac{2}{\eta} \|v_0\|_{L^2}^2 + \frac{(1{+}\eta)\alpha^2}{4 
  \pi} \log (N t)~.
\end{align*}
Similarly, using \eqref{uchi5}, we find
\begin{align*}
 \int_0^t \|\nabla v(s)\|_{L^2}^2 \dd s  \,&\le\, 2 \int_0^t 
  \|\nabla \tilde v(s)\|_{L^2}^2 \dd s  + 2\alpha^2 \int_0^t 
  \|\nabla u^\chi(s)- \nabla u^\chi(s+\tau)\|_{L^2}^2 \dd s \\
  \,&\le\, 2 \int_0^t \|\nabla \tilde v(s)\|_{L^2}^2 \dd s 
  +  \frac{\alpha^2}{4\pi} \,\log (1 + t) + C\rho^2 \alpha^2~.
\end{align*}
Thus, it follows from \eqref{tildevpoly} that 
\begin{align}\label{logprelim1}
  \|v(t)\|_{L^2}^2 \,&\le\, \frac{(1{+}\eta)^3\alpha^2}{4\pi} \log t
  + \frac{C}{\eta}\Bigl(\|v_0\|_{L^2}^2 + \alpha^2 \log (N+1) + \alpha^2
  \rho^2\Bigr)~, \\ \label{logprelim2} 
  \int_0^t \|\nabla v(s)\|_{L^2}^2 \dd s  \,&\le\, 
  \frac{(1{+}\eta)^3\alpha^2}{2\pi} \log(1+t ) + \frac{C}{\eta} 
  \Bigl(\|v_0\|_{L^2}^2 + \alpha^2 \rho^2\Bigr) + 
  C\alpha^2 \log N~,
\end{align}
for some universal constant $C > 0$. Setting $\eta = 1$ and 
using the definition of $N$, we see that \eqref{vlog} follows 
from \eqref{logprelim1}, \eqref{logprelim2}. \QED

\section{Estimate for a fractional primitive of the velocity 
field}\label{sec5}

In this final section, we consider the solution of \eqref{veq2} with
initial data $v_0 \in L^2_\sigma(\Omega) \cap L^q(\Omega)^2$, for some
fixed $q \in (1,2)$, and we denote $\mu = 1/q - 1/2 \in (0,1/2)$. If
$A$ is the Stokes operator in $L^2_\sigma(\Omega)$, we recall that $A$
is selfadjoint and nonnegative in $L^2_\sigma(\Omega)$, so that the
fractional power $A^\beta$ can be defined for all $\beta > 0$. The
following result shows that the range of $A^\mu$ contains the (dense)
subspace $L^2_\sigma(\Omega) \cap L^q(\Omega)^2$.

\begin{lem}\label{Amu} {\rm \cite{BM,KO2}} Let $q \in (1,2)$ and
$\mu = 1/q - 1/2$. For all $v \in L^2_\sigma(\Omega) \cap L^q(\Omega)^2$, 
there exists a unique $w \in D(A^\mu) \subset L^2_\sigma(\Omega)$ such 
that $v = A^\mu w$. Moreover, there exists a constant $C=C(q) > 0$ 
(independent of $v$ and $\Omega$) such that $\|w\|_{L^2(\Omega)} \le 
C\|v\|_{L^q(\Omega)}$. 
\end{lem}

\begin{rem}\label{A-mu}
If $v,w$ are as in Lemma~\ref{Amu}, we denote $w = A^{-\mu}v$. 
The fact that inequality $\|w\|_{L^2(\Omega)} \le C\|v\|_{L^q(\Omega)}$ 
holds with a constant $C$ independent of the domain $\Omega$ follows 
directly from the proof given in \cite[Lemmas 2.1 and 2.2]{KO2}. 
\end{rem}

As a first application of Lemma~\ref{Amu}, we give a short proof
of inequality \eqref{Stokesbd}, which was used in Section~\ref{sec3}. 

\medskip\noindent{\bf Proof of Proposition~\ref{Stokes}.} It is 
sufficient to prove \eqref{Stokesbd} for $1 < q < 2$. Let $\mu = 
1/q-1/2$, and let $v_0\in L^2_\sigma(\Omega) \cap L^q(\Omega)^2$. 
By Lemma~\ref{Amu}, there exists a unique $w_0\in D(A^\mu)$ such that 
$v_0 = A^\mu w_0$. Thus
\begin{align*}
  \|S(t)v_0\|_{L^2(\Omega )} \,=\, \|A^\mu S(t)w_0\|_{L^2 (\Omega )} \,\le\,  
  t^{-\mu} \|w_0\|_{L^2(\Omega)} \,\le\, C t^{-\mu} \|v_0\|_{L^q(\Omega)}~,
\end{align*}
with $C$ depending only on $q$. The estimate for the first derivative
is proved in the same way, since $\|\nabla S(t)v_0\|_{L^2(\Omega)} = 
\|A^{\mu+1/2} S(t)w_0\|_{L^2(\Omega )}$. This proves \eqref{Stokesbd} for 
all $v_0\in L^2_\sigma(\Omega) \cap L^q(\Omega)^2$, and the general case
follows by a density argument. \QED

\medskip 
Let $v \in C^0([0,\infty); L^2_\sigma(\Omega)) \cap C^0((0,\infty);
H^1_0(\Omega)^2)$ be the solution of \eqref{veq2} with initial data
$v_0$, which was constructed in Sections~\ref{sec3} and
\ref{sec4}. Since $v_0 \in L^q_\sigma(\Omega)$ by assumption, it is
rather straightforward to verify that $v(t) \in L^q_\sigma(\Omega)$
for all $t > 0$. Thus, by Lemma~\ref{Amu}, we can define $w(t) =
A^{-\mu} v(t)$ for all $t > 0$. This quantity solves the equation
\begin{equation}\label{weq}
  \partial_t w + A w + \alpha F_\mu(u^\chi,v) + \alpha F_\mu(v,u^\chi) +
  F_\mu(v,v) \,=\, \alpha A^{-\mu}R^\chi~,
\end{equation}
where $F_\mu(u,v)$ is the bilinear term formally defined by
\begin{equation}\label{Fmudef}
  F_\mu(u,v) \,=\, A^{-\mu} P(u\cdot\nabla)v~.
\end{equation}
We refer to \cite[Section 2]{KO2} for a rigorous definition and 
a list of properties of the bilinear map $F_\mu$.  Our goal here
is to establish the following estimate\:

\begin{prop}\label{wenergy}
There exist positive constants $K_2$ and $c$ (independent of 
$\Omega$) such that, for any $\alpha \in \R$ and any solution 
$v$ of \eqref{veq2} with initial data $v_0 \in L^2_\sigma(\Omega) 
\cap L^q(\Omega)^2$, the function $w(t) = A^{-\mu}v(t)$ satisfies, 
for all $t \ge 1$, 
\begin{equation}\label{wbound}
  \|w(t)\|_{L^2}^2 + \int_0^t \|\nabla w(s)\|_{L^2}^2\dd s \,\le\, 
  K_2 (1+t)^{c \alpha^2} \!\exp\Bigl(K_2(\|v_0\|_{L^2}^2 + 
  D_{\alpha,\rho})\Bigr)(\|v_0\|_{L^q}^2 + \rho^2 \alpha^2)~,
\end{equation}
where $D_{\alpha,\rho} = \alpha^2\log(1+|\alpha|) + \alpha^2 \rho^2$.
\end{prop}

\noindent{\bf Proof.} Taking the scalar product of both sides
of \eqref{weq} with $w$, we obtain
\begin{align}\nonumber
  \frac12 \frac{\D}{\D t} \|w(t)\|_{L^2}^2 + \|A^{1/2} w(t)\|_{L^2}^2 
  &+ \alpha \langle F_\mu(u^\chi(t),v(t)),w(t)\rangle +  \alpha
  \langle F_\mu(v(t),u^\chi(t)),w(t)\rangle \\ \label{wfirst} 
  &+ \langle F_\mu(v(t),v(t)),w(t)\rangle 
  \,=\, \alpha \langle A^{-\mu}R^\chi(t),w(t)\rangle~.
\end{align}
We recall that $\|A^{1/2} w\|_{L^2} = \|\nabla w\|_{L^2}$ for all $w \in 
D(A^{1/2}) = L^2_\sigma(\Omega) \cap H^1_0(\Omega)^2$. To bound the 
other terms, we observe that
\begin{align*}
   |\langle F_\mu(u^\chi,v),w\rangle| \,&=\, |\langle (u^\chi\cdot \nabla) 
   v,A^{-\mu}w\rangle|  \,=\, |\langle (u^\chi\cdot \nabla)A^{-\mu}w 
   ,v\rangle| \\ \,&\le\, \|u^\chi\|_{L^\infty} \|A^{\frac12-\mu}w\|_{L^2} 
  \|v\|_{L^2} \,=\,  \|u^\chi\|_{L^\infty} \|A^{\frac12-\mu}w\|_{L^2} 
  \|A^{\mu}w\|_{L^2} \\ 
  \,&\le\, \|u^\chi\|_{L^\infty} \|A^{1/2}w\|_{L^2} \|w\|_{L^2}~,
\end{align*}
where in the last inequality we used the interpolation inequality 
for fractional powers of $A$. The same argument shows that $|\langle F_\mu
(v,u^\chi),w\rangle| \le \|u^\chi\|_{L^\infty} \|A^{1/2}w\|_{L^2} \|w\|_{L^2}$. 
In a similar way, we find
\begin{align*}
  |\langle F_\mu(v,v),w\rangle| \,&=\, |\langle (v\cdot \nabla) 
   v,A^{-\mu}w\rangle|  \,=\, |\langle (v\cdot \nabla)A^{-\mu}w 
   ,v\rangle| \\ \,&\le\, \|v\|_{L^4}^2 \|A^{\frac12-\mu}w\|_{L^2} 
   \,\le\, C_*^2 \|\nabla v\|_{L^2} \|v\|_{L^2} \|A^{\frac12-\mu}w\|_{L^2} 
   \\ \,&\le\, C_*^2 \|\nabla v\|_{L^2} \|A^{1/2}w\|_{L^2} \|w\|_{L^2}~,
\end{align*}
where $C_* > 0$ is the best constant of Gagliardo-Nirenberg's inequality
\begin{equation}\label{gagliardo}
  \|f\|_{L^4(\R^2) } \,\le\, C_* \|f\|_{L^2 (\R^2 )}^{\frac12} \|\nabla 
  f\|_{L^2(\R^2) }^{\frac12}~.
\end{equation}
Finally, since $|\langle A^{-\mu}R^\chi,w \rangle| =  |\langle R^\chi,
A^{-\mu} w \rangle| \le \kappa_2\rho (1+t)^{-1}\|A^{\frac12-\mu}w\|_{L^2}$ 
by \eqref{Rchi2}, we can use interpolation and Young's inequality
to obtain
\[
  |\alpha\langle A^{-\mu}R^\chi,w \rangle|  \,\le\, \frac{\kappa_2 \rho 
  |\alpha|}{1+t}\|A^{1/2}w\|_{L^2}^{1-2\mu}\|w\|_{L^2}^{2\mu} \,\le\, 
  \frac\eta4 \|A^{1/2}w\|_{L^2}^2 + \frac{\|w\|_{L^2}^2}{2(1+t)^{\gamma_1}} + 
  \frac{C_\eta\rho^2\alpha^2}{2(1+t)^{\gamma_2}}~,
\]
for some exponents $\gamma_1,\gamma_2 > 1$ satisfying $\gamma_2 + 
2\mu\gamma_1 = 2$. Here $\eta \in (0,1]$ is as in the proof of 
Proposition~\ref{logbound}, and $C_\eta > 0$ denotes a constant 
depending only on $\eta$. Inserting all these estimates into 
\eqref{wfirst}, we arrive at
\begin{equation}\label{wsecond}
  \frac{\D}{\D t} \|w\|_{L^2}^2 +  2 \|\nabla w\|_{L^2}^2 \,\le\, 
  2 H \|\nabla w\|_{L^2} \|w\|_{L^2} + \frac{\eta}{2}\|\nabla w\|_{L^2}^2
  + \frac{\|w\|_{L^2}^2}{(1+t)^{\gamma_1}} + \frac{C_\eta \rho^2 \alpha^2}
  {(1+t)^{\gamma_2}}~,  
\end{equation}
where $H = 2 |\alpha| \|u^\chi\|_{L^\infty} + C_*^2 \|\nabla v\|_{L^2}$.

To exploit \eqref{wsecond}, we apply Young's inequality again and 
obtain the differential inequality
\[
  \frac{\D}{\D t} \|w\|_{L^2}^2 +  \eta \|\nabla w\|_{L^2}^2 \,\le\, 
  \Bigl(\frac{H^2}{2-3\eta/2} + \frac{1}{(1+t)^{\gamma_1}}\Bigr)
  \|w\|_{L^2}^2 + \frac{C_\eta \rho^2 \alpha^2}{(1+t)^{\gamma_2}}~,  
\]
which can be integrated using Gronwall's lemma. The result is
\begin{equation}\label{wthird}
  \|w(t)\|_{L^2}^2 + \eta \int_0^t \|\nabla w(s)\|_{L^2}^2\dd s 
  \,\le\, C \exp\Bigl(\frac{\Phi(t)}{1-3\eta/4}\Bigr) 
  \Bigl(\|w_0\|_{L^2}^2 + C_\eta \rho^2 \alpha^2\Bigr)~, \quad t \ge 0~,
\end{equation}
where $\Phi(t) = \frac12 \int_0^t H(s)^2\dd s$ and $C$ is a positive 
constant depending only on $\gamma_1, \gamma_2$. It remains to estimate
the quantity $\Phi(t)$ in \eqref{wthird}. Using \eqref{uchi1} with 
$p=\infty$, the logarithmic energy estimate \eqref{logprelim2}, 
and Minkowski's inequality, we find
\begin{align}\nonumber
  2 \Phi(t) \,&=\, \int_0^t H(s)^2 \dd s \,\le\, \int_0^t \Bigl\{
  \frac{2|\alpha|a_\infty} {(1+s)^{1/2}} + C_*^2 \|\nabla v(s)\|_{L^2} 
  \Bigr\}^2 \dd s \\ \label{Phiest}
  \,&\le\, \Bigl\{|\alpha| \log(1+t)^{1/2} \Bigl(2 a_\infty + \frac{C_*^2
  (1+\eta)^{\frac32}}{\sqrt{2 \pi}}\Bigr) + C_\eta (\|v_0\|_{L^2} + 
  D_{\alpha,\rho}^{1/2})\Bigr\}^2 \\ \nonumber
   \,&\le\, 2 C_0 (1+\eta)^4 \alpha^2 \log(1+t) +  C_\eta (\|v_0\|_{L^2}^2 + 
  D_{\alpha,\rho})~, \qquad t \ge 1~,
\end{align}
where $D_{\alpha,\rho} = \alpha^2\log(1+|\alpha|) + \alpha^2 \rho^2$ and 
\begin{equation}\label{C0def}
  C_0 \,=\, \frac12 \Bigl(2 a_\infty + \frac{C_*^2}{\sqrt{2\pi}}\Bigr)^2~.
\end{equation}
If we now replace \eqref{Phiest} into \eqref{wthird} and set $\eta = 1$, 
we obtain \eqref{wbound} since $\|w_0\|_{L^2} \le C\|v_0\|_{L^q}$ by 
Lemma~\ref{Amu}. This concludes the proof. \QED

\begin{cor}\label{vcor}
Under the assumptions of Proposition~\ref{wenergy}, there exists a 
positive constant $K$ depending on $\Omega$, $\alpha$, $q$, and 
$\|v_0\|_{L^2\cap L^q}$ such that, for any $T \ge 2$, there exists a 
time $t \in [T/2,T]$ for which
\begin{equation}\label{vest}
  \|v(t)\|_{L^2(\Omega)}^2 \,\le\, K (1+t)^{c\alpha^2 - 2\mu}~.
\end{equation}
\end{cor}

\noindent{\bf Proof.} Fix $T \ge 2$. In view of \eqref{wbound}, there 
exists a time $t \in [T/2,T]$ such that
\begin{align*}
  \|\nabla w(t)\|_{L^2}^2 \,\le\, \frac{2}{T}\int_{T/2}^T 
  \|\nabla w(s)\|_{L^2}^2\dd s \,\le\, \frac{2}{T}\,C (1+T)^{c\alpha^2} 
  \,\le\, 2^{c\alpha^2+2} C(1+t)^{c\alpha^2-1}~,
\end{align*}
where $C$ depends on $\rho$, $\alpha$, $q$ and $\|v_0\|_{L^2 \cap L^q}$. 
Moreover, $\|w(t)\|_{L^2}^2 \le C(1+t)^{c\alpha^2}$ by \eqref{wbound}.
Thus, using the interpolation inequality $\|v(t)\|_{L^2} = \|A^\mu 
w(t)\|_{L^2} \le \|\nabla w(t)\|_{L^2}^{2\mu} \,\|w(t)\|_{L^2}^{1-2\mu}$, 
we obtain \eqref{vest}.
\QED

\bigskip\noindent{\bf Proof of Theorem~\ref{main}.}  Fix $q \in 
(1,2)$, and assume that $\epsilon > 0$ is small enough so that
$c\epsilon^2 < 2\mu$, where $\mu = 1/q-1/2$ and $c$ is as in
Proposition~\ref{wenergy}. We also suppose that $\epsilon \le
\delta/2$, where $\delta > 0$ is as in Proposition~\ref{small}.  Given
$\alpha \in [-\epsilon,\epsilon]$ and $v_0 \in L^2_\sigma(\Omega) \cap
L^q(\Omega)^2$, let $v \in C^0([0,\infty); L^2_\sigma(\Omega)) \cap
C^0((0,\infty); H^1_0(\Omega)^2)$ be the solution of \eqref{veq2} with
initial data $v(0) = v_0$, which was constructed in
Sections~\ref{sec3} and \ref{sec4}. In view of \eqref{vest}, since 
$c\alpha^2 < 2\mu$, we can take $t_0 > 0$ large enough (depending
on $\Omega$, $\alpha$, and $v_0$) so that $\|v(t_0)\|_{L^2} \le
V_\Omega$, where $V_\Omega$ is as in Proposition \ref{small}. Moreover, 
since $v(t_0) = A^\mu w(t_0)$ for some $w(t_0) \in L^2_\sigma(\Omega)$, 
we have
\[
  \sup_{\tau > 0} \tau^\mu\|S(\tau)v_0\|_{L^2} +  \sup_{\tau > 0}
  \tau^{\mu + \frac12} \|\nabla S(\tau)v_0\|_{L^2} \,\le\, 
  C \|w(t_0)\|_{L^2} \,<\, \infty~.
\]
Applying Proposition~\ref{small}, we conclude that the solution 
$v$ of \eqref{veq2} satisfies \eqref{vdecay}, namely
\begin{equation}\label{conv2}
  \|u(\cdot,t) -\alpha u^\chi(\cdot,t)\|_{L^2} + 
  t^{1/2}\|\nabla u(\cdot,t) -\alpha\nabla u^\chi(\cdot,t)
  \|_{L^2(\Omega)} \,=\, \OO(t^{-\mu})~,
\end{equation}
as $t \to \infty$. But $\|u^\chi - \Theta\|_{L^2} + \|\nabla u^\chi
- \nabla \Theta\|_{L^2} \le C(1+t)^{-1}$ for all $t \ge 0$, hence 
\eqref{conv} follows from \eqref{conv2}. 
\QED

\bigskip\noindent{\bf Proof of Corollary~\ref{explicit}.}  
The proof of Proposition~\ref{wenergy} shows that the constant
$c$ in \eqref{wbound}, \eqref{vest} satisfies $c \le C_0(1+\OO(\eta))$, 
where $C_0$ is defined in \eqref{C0def} and $\eta \in (0,1]$ can
be chosen arbitrarily small. On the other hand, since by assumption 
$v_0 \in L^2_\sigma(\Omega) \cap L^q(\Omega)^2$ for all $q \in (1,2)$, we 
can take $\mu = 1/q - 1/2$ arbitrarily close to $1/2$. Thus, if we 
assume that $|\alpha| < \epsilon_* = C_0^{-1/2}$, we see that the 
condition $c\alpha^2 < 2\mu$ can be fulfilled by an appropriate 
choice of $\eta$ and $\mu$. Now, take $t \ge 2$ and let $t_0 \in [t/2,t]$ 
be the time defined in Corollary \ref{vcor}, for which $\|v(t_0)\|_{L^2}^2 
\le K (1+t_0)^{c\alpha^2 - 2\mu}$. Using \eqref{vpoly} with $\eta = 1$, 
we conclude
\[
  \|v(t)\|_{L^2}^2 \,\le\, C\Bigl(\frac{1+t}{1+t_0}\Bigr)^{2 b_\infty 
  |\alpha|} \Bigl(\|v(t_0)\|_{L^2}^2 + (1+t_0)^{-1}\Bigr) \,\le\, C 
  (1+t)^{c\alpha^2 - 2\mu} \,\xrightarrow[t \to \infty]{}\, 0~,
\]
which is the desired result. Here the constant $C > 0$ depends 
on $\alpha$, $\rho$, and $v_0$, but not on $t$. To estimate
$\epsilon_*$, we use \eqref{C0def} and observe that $a_\infty = 
\|\Theta_0\|_{L^\infty} \approx 0.050784$. Morevoer, the optimal constant 
in the Gagliardo-Nirenberg inequality \eqref{gagliardo} satisfies 
$C_*^4 \le 2/(3\pi)$, see \cite{DD}. Using these values, we find 
$C_0 \le 0.0407108$, hence $\epsilon_* = C_0^{-1/2} \ge 4.95616$. 
Finally, it was kindly pointed to us by Jean Dolbeault that the 
optimal constant $C_*$ can be computed numerically\: $C_* \approx 
0.6430$. This yields the approximate value $\epsilon_* \approx 5.306$. 
\QED

\section{Appendix\: Proof of Proposition~\ref{vorticity}}

We recall the following characterization of the space $\dot W_{0,\sigma}^{1,p}
(\Omega)$ for $1 \le p < 2$\:
\begin{align}\label{hatW}
  \dot W_{0,\sigma}^{1,p}(\Omega) \,=\, \big \{ u \in L^\frac{2 p}{2-p}
  (\Omega )^2 ~\big |~ \|\nabla u\|_{L^p}<\infty\,,~ u=0~\hbox{on}~
  \partial\Omega\,,~\div u = 0 ~\hbox{in}~\Omega \big \}~,
\end{align}
see e.g. \cite[Chapter III.5]{Galdi}. Here $\nabla u$ and $\div u$ denote 
weak derivatives of $u$, and the condition ``$u=0$ on $\partial \Omega$''
means that the boundary trace of $u$, which is well defined because
$\nabla u \in L^p(\Omega)^4$, vanishes.

Given $u_0 \in \dot W_{0,\sigma}^{1,p}(\Omega)$ satisfying \eqref{omloc}, 
we define $u : \R^2 \to \R^2$ and $\omega : \R \to \R^2$ by
\[
  u(x) \,=\, \left\{ \begin{array}{ccc} \!\!u_0(x) & \hbox{if} & x \in 
  \Omega~, \\ 0 & \hbox{if} & x \notin \Omega~,\end{array}
  \right. \qquad
  \omega(x) \,=\, \left\{ \begin{array}{ccc} \!\!\omega_0(x) & 
  \hbox{if} & x \in \Omega~, \\ 0 & \hbox{if} & x \notin \Omega~.
  \end{array}\right.
\]
Since $u = 0$ on $\partial \Omega$, we have $\nabla u \in L^p(\R^2)^4$ and 
$\partial_1 u_2 - \partial_2 u_1 = \omega \in L^p(\R^2)$. Moreover 
\eqref{omloc} implies that $\omega \in L^2(m)$ for some $m > 2/q > 1$, 
where 
\[
  L^2(m) \,=\, \Bigl\{\omega \in L^2(\R^2) \,\Big|\,  \int_{\R^2} 
  (1+|x|^2)^m |\omega(x)|^2 \dd x \,<\, \infty \Bigr\}~. 
\]
Thus, using H\"older's inequality, it is easy to verify that $\omega \in 
L^1(\R^2)$, so that we can define
\[
  \alpha \,=\, \int_{\R^2} \omega(x) \dd x \,=\, \int_{\Omega} 
  \omega_0(x) \dd x~.
\]
Moreover, using the the Biot-Savart formula in $\R^2$ and the fact that 
$u \in L^{2p/(p-2)}(\R^2)^2$, we obtain the equality
\begin{equation}\label{BS}
  u(x) \,=\, \frac{1}{2\pi} \int_{\R^2} \frac{(x-y)^\perp}{|x-y|^2}
  \,\omega(y)\dd y \,=\, \int_{\Omega} \frac{(x-y)^\perp}{|x-y|^2}
  \,\omega_0(y)\dd y~, 
\end{equation}
for almost all $x \in \R^2$.  We emphasize at this point that the
representation \eqref{BS} is {\em not} what is usually called the
Biot-Savart law in the domain $\Omega$, because the velocity field
defined by \eqref{BS} for an arbitrary vorticity $\omega_0 \in
L^1(\Omega)$ will not, in general, be tangent to the boundary on
$\partial \Omega$.  However, if we start from a velocity field $u_0$
that vanishes on $\partial \Omega$, the argument above shows that
\eqref{BS} holds with $\omega_0 = \curl u_0$. We refer to \cite{ILN}
for a more detailed discussion of the Biot-Savart law in a
two-dimensional exterior domain.

Now, we decompose
\[
  u(x) \,=\, \alpha u^\chi(x,0) + v(x)~, \qquad
  \omega(x) \,=\, \alpha \omega^\chi(x,0) + w(x)~, 
  \qquad x \in \R^2~,  
\]
where $u^\chi$, $\omega^\chi$ are defined in \eqref{uchidef}, 
\eqref{omegachi}. By construction, we have $w \in L^2(m)$ for
some $m \in (1,2)$ and $\int_{\R^2} w\dd x = 0$. Applying 
\cite[Proposition~B.1]{GW1}, we deduce that the corresponding 
velocity field $v$, which is obtained from $w$ via the Biot-Savart 
law in $\R^2$, satisfies
\[
  \int_{\R^2} (1+|x|^2)^{\frac{mr}{2}-1}|v(x)|^r \dd x \,<\, \infty~,
\]
for all $r > 2$. Using H\"older's inequality again, we conclude
that $v \in L^s(\R^2)^2$ for all $s > 2/m$, hence in particular 
$v \in  L^2(\R^2)^2 \cap L^q(\R^2)^2$. Clearly $v(x) = 0$ for all
$x \notin \Omega$, hence denoting by $v_0$ the restriction 
of $v$ to $\Omega$ we obtain \eqref{init} with $v_0 \in 
L^2_\sigma(\Omega) \cap L^q(\Omega)^2$. \QED

\bigskip\noindent{\bf Acknowledgements.} The authors are 
partially supported by the ANR project PREFERED (Th.G.) and
by the Grant-in-Aid for Young Scientists (B) 22740090 (Y.M.).
This work was initiated during a visit of Y.M. at Universit\'e 
Joseph Fourier (Grenoble I), whose hospitality is gratefully 
acknowledged.

\end{document}